\documentclass[11pt]{amsart}  
\usepackage{amsmath,amssymb,amsfonts,times,hyperref}

\title{Symmetry in semidefinite programs}

\author{Frank Vallentin} 

\address{F. Vallentin, Centrum voor Wiskunde en Informatica (CWI),
Kruislaan 413, 1098 SJ Amsterdam, The Netherlands}

\email{f.vallentin@cwi.nl}

\date{July 30, 2008}

\thanks{The author was supported by the Netherlands Organization for
Scientific Research under grant NWO 639.032.203 and by the Deutsche
Forschungsgemeinschaft (DFG) under grant SCHU 1503/4-2.}


\date{\today}
\newtheorem{defi}{Definition}[section]

\newtheorem{proposition}[defi]{Proposition}
\newtheorem{theorem}[defi]{Theorem}
\newtheorem{remark}[defi]{Remark}

\newtheorem{lemma}[defi]{Lemma}

\newcommand{\R}{{\mathbb{R}}} 
\newcommand{\C}{{\mathbb{C}}}

\newcommand{\MA}{{\mathcal{A}}} 
\newcommand{\MB}{{\mathcal{B}}}

\newcommand{\Ort}{{\operatorname{O}}}
\newcommand{\Uni}{{\operatorname{U}}}

\newcommand{\trace}{\operatorname{trace}}

\newcommand{\Comm}{\operatorname{Comm}}

\newcommand{\prodtrace}[2]{\langle #1, #2 \rangle}

\begin{document}

\begin{abstract}
  This paper is a tutorial in a general and explicit procedure to
  simplify semidefinite programs which are invariant under the action
  of a symmetry group. The procedure is based on basic notions of
  representation theory of finite groups. As an example we derive the
  block diagonalization of the Terwilliger algebra of the binary
  Hamming scheme in this framework.  Here its connection to the
  orthogonal Hahn and Krawtchouk polynomials becomes visible.
\end{abstract}

\maketitle

\section{Introduction}

A \textit{(complex) semidefinite program} is an optimization problem
of the form
\begin{equation}
\label{sdp standard}
\max\{\prodtrace{C}{Y} : \text{$\prodtrace{A_i}{Y} = b_i$, $i =
1, \ldots, n$, and $Y \succeq 0$}\},
\end{equation}
where $A_i \in \C^{X \times X}$, and $C \in \C^{X \times X}$ are given
Hermitian matrices whose rows and columns are indexed by a finite set
$X$, $(b_1, \ldots, b_n)^t \in \R^n$ is a given vector and $Y \in
\C^{X \times X}$ is a variable Hermitian matrix and where ``$Y \succeq
0$'' means that $Y$ is positive semidefinite. Here $\prodtrace{C}{Y} =
\trace(CY)$ denotes the trace product between symmetric matrices.

Semidefinite programming is an extension of linear programming and has
a wide range of applications: combinatorial optimization and control
theory are the most famous ones. Although semidefinite programming has
an enormous expressive power in formulating convex optimization
problems it has a few practical drawbacks: Highly robust and highly
efficient solvers, unlike their counterparts for solving linear
programs, are currently not available. So it is crucial to exploit the
problems' structure to be able to perform computations.

In the last years many results were obtained if the problem under
consideration has symmetry. This was done for a variety of problems
and applications: interior point algorithms (Kanno, Ohsaki, Murota,
Katoh \cite{komk-2001} and de Klerk, Pasechnik \cite{dp-2005}),
polynomial optimization (Parrilo, Gatermann \cite{gp-2004} and
Jansson, Lasserre, Riener, Theobald \cite{jlrt-2006}), truss topology
optimization (Bai, de Klerk, Pasechnik, Sotirov \cite{bdps-2007}),
quadratic assignment (de Klerk, Sotirov \cite{ds-2007}), fast mixing
Markov chains on graphs (Boyd, Diaconis, Xiao \cite{bdx-2004}), graph
coloring (Gvozdenovi\'c, Laurent \cite{gl-2007}), crossing numbers for
complete binary graphs (de Klerk, Pasechnik, Schrijver
\cite{dps-2007}) and coding theory (Schrijver \cite{schrijver-2005},
Gijswijt, Schrijver, Tanaka \cite{gst-2006} and Laurent
\cite{laurent-2007}).

In all these applications the underlying principles are similar: one
simplifies the original semidefinite program which is invariant under
a group action by applying an algebra isomorphism mapping a ``large''
matrix algebra to a ``small'' matrix algebra. Then it is sufficient to
solve the semidefinite program using the smaller matrices. The
existence of an appropriate algebra isomorphism is a classical fact
from Artin-Wedderburn theory. However, in the above mentioned papers
the explicit determination of an appropriate isomorphism is rather
mysterious. The aim of this paper is to give an algorithmic way to do
this which also is well-suited for symbolic calculations by hand.

The paper is structured as follows: Section~\ref{background} recalls
basic definitions and shows how the Artin-Wedderburn theorem stated in
\eqref{isomorphism} can be applied to simplify a semidefinite program
invariant under a group action. In Section~\ref{blockdiag} we construct an
explicit algebra isomorphism. In Section~\ref{terwilliger} we apply
this to the Terwilliger algebra of the binary Hamming scheme.

This paper is of expository nature and probably few of the results are
new. On the other hand a tutorial of how to use symmetry in
semidefinite programming is not readily available. Furthermore our
treatment of the Terwilliger algebra for binary codes provides an
alternative point of view which emphasizes the action of the symmetric
group. Schrijver \cite{schrijver-2005} treated the Terwilliger algebra
with elementary combinatorial and linear algebraic arguments. Our
derivation has the advantage that it gives an interpretation for the
matrix entries in terms of Hahn polynomials. In a similar way one can
derive the block diagonalization of the Terwilliger algebra for
nonbinary codes which was computed by Gijswijt, Schrijver, Tanaka
\cite{gst-2006}. Here products of Hahn and Krawtchouk polynomials
occur.

\section{Background and notation}
\label{background}

In this section we present the basic framework for simplifying a
semidefinite program invariant under a group action.

Let $G$ be a finite group which acts on a finite set $X$ by $(a,x)
\mapsto ax$ with $a \in G$ and $x \in X$. This group action extends to
an action on pairs $(x, y) \in X \times X$ by $(a, (x,y)) \mapsto (ax,
ay)$. In this way it extends to square matrices whose rows and columns
are indexed by $X$: for an $X\times X$-matrix $M$ we have $aM(x,y) =
M(ax,ay)$. Here $M(x,y)$ denotes the entry of $M$ at position
$(x,y)$. A matrix $M$ is called \textit{invariant under $G$} if $M =
aM$ for all $a \in G$.

A Hermitian matrix $Y \in \C^{X \times X}$ is called a
\textit{feasible solution} of \eqref{sdp standard} if it fulfills the
conditions $\prodtrace{A_i}{Y} = b_i$ and $Y \succeq 0$. It is called
an \textit{optimal solution} if it is feasible and if for all other
feasible solutions $Y'$ we have $\prodtrace{C}{Y} \geq
\prodtrace{C}{Y'}$. In the following we assume that the semidefinite
program \eqref{sdp standard} has an optimal solution.

We say that the semidefinite program~\eqref{sdp standard} is
\textit{invariant under $G$} if for every feasible solution $Y$ and
for every $a \in G$ the matrix $aY$ is again a feasible solution and
if it is satisfies $\prodtrace{C}{aY} = \prodtrace{C}{Y}$ for all $a
\in G$. Because of the convexity of \eqref{sdp standard}, one can find
an optimal solution of \eqref{sdp standard} in the subspace $\MB$ of
matrices which are invariant under $G$. In fact, if $Y$ is an optimal
solution of \eqref{sdp standard}, so is its \textit{group average}
$\frac{1}{|G|} \sum_{a \in G} aY$. Hence, \eqref{sdp standard} is
equivalent to
\begin{equation}
\label{sdp standard invariant}
\max\{\prodtrace{C}{Y} : \mbox{$\prodtrace{A_i}{Y} = b_i$, $i =
1, \ldots, n$, $Y \succeq 0$, and $Y \in \MB$}\}.
\end{equation}

The set $X \times X$ can be decomposed into the orbits $R_1, \ldots,
R_N$ by the action of $G$. For every $r \in \{1, \ldots, N\}$ we
define the matrix $B_r \in \{0,1\}^{X \times X}$ by $B_r(x, y) = 1$ if
$(x,y) \in R_r$ and $B_r(x,y) = 0$ otherwise. Then $B_1, \ldots, B_N$
forms a basis of $\MB$. We call $B_1, \ldots, B_N$ the
\textit{canonical basis} of $\MB$. If $(x,y) \in R_r$ we also write
$B_{[x,y]}$ instead of $B_r$. Note that $B_{[y,x]}$ is the transpose
of the matrix $B_{[x,y]}$.

So the first step to simplify a semidefinite program which
is invariant under a group is as follows:

\textit{\noindent
If the semidefinite program~\eqref{sdp standard} is
invariant under $G$, then~\eqref{sdp standard} is equivalent to
\begin{equation}
\label{sdp first reduction}
\begin{array}{lcl}
\max\Big\{c_1 y_1 + \cdots + c_N y_N & : & y_1, \ldots, y_N \in \C,\\
& & a_{i1}y_1 + \cdots + a_{iN}y_N = b_i,\; i = 1, \ldots, n,\\
& & \mbox{$y_j = \overline{y_{k}}$ if $B_j = (B_{k})^t$,}\\
& & y_1 B_1 + \cdots + y_N B_N \succeq 0\Big\},\\
\end{array}
\end{equation}
where  $c_r = \prodtrace{C}{B_r}$, and $a_{ir} = \prodtrace{A_i}{B_r}$.
}

The following obvious property is crucial for the next step of
simplifying \eqref{sdp first reduction}: The subspace $\MB$ is closed
under matrix multiplication. So $\MB$ is a (semisimple) algebra over
the complex numbers. The Artin-Wedderburn theory (cf. \cite[Chapter
1]{lam-1991}) gives:

\textit{
\noindent
There are numbers $d$, and $m_1, \ldots, m_d$ so that there is an
algebra isomorphism
\begin{equation}
\label{isomorphism}
\varphi : \MB \to \bigoplus_{k=1}^d \C^{m_k \times m_k}.
\end{equation}
}

This applied to \eqref{sdp first reduction} gives the final step of
simplifying \eqref{sdp standard}:

\textit{
\noindent
If the semidefinite program~\eqref{sdp standard} is
invariant under $G$, then~\eqref{sdp standard} is equivalent to
\begin{equation}
\label{sdp block}
\begin{array}{lcl}
\max\Big\{c_1 y_1 + \cdots + c_N y_N & : & y_1, \ldots, y_N \in \C,\\
& & a_{i1}y_1 + \cdots + a_{iN}y_N = b_i,\; i = 1, \ldots, n,\\
& & \mbox{$y_j = \overline{y_{k}}$ if $B_j = (B_{k})^t$,}\\
& & y_1 \varphi(B_1) + \cdots + y_N \varphi(B_N) \succeq 0\Big\}.\\
\end{array}
\end{equation}
}

Notice that since $\varphi$ is an algebra isomorphism between matrix
algebras with unity, $\varphi$ preserves eigenvalues and hence
positive semidefiniteness. In accordance to the literature, applying
$\varphi$ to a semidefinite program is called
\textit{block diagonalization}.

The advantage of \eqref{sdp block} is that instead of dealing with
matrices of size $|X| \times |X|$ one has to deal with block diagonal
matrices with $d$ block matrices of size $m_1, \ldots, m_d$,
respectively. In many applications the sum $m_1 + \cdots + m_d$ is
much smaller than $|X|$ and in particular many practical solvers take
advantage of the block structure to speed up the numerical
calculations.

\section{Determining a block diagonalization}
\label{blockdiag}

In this section we give an explicit construction of an algebra
isomorphism $\varphi$. It has two main features: One can turn the
construction into an algorithm as we show at the end of this section,
and one can use it for symbolic calculations by hand as we demonstrate
in Section~\ref{terwilliger}.

\subsection{Construction}

We begin with some basic notions from representation theory of finite
groups.  Consider the complex vector space $\C^X$ of vectors indexed
by $X$ with inner product $(f, g) = \frac{1}{|X|} \sum_{x \in X} f(x)
\overline{g(x)}$. The group $G$ acts on $\C^X$ by $af(x) =
f(a^{-1}x)$. Note that the inner product on $\C^X$ is invariant under
the group action: For all $f, g \in \C^X$ and all $a \in G$ we have
$(af, ag) = (f, g)$. A subspace $H \subseteq \C^X$ is called a
\textit{$G$-space} if $GH \subseteq H$ where $GH = \{af : f \in H, a
\in G\}$. It is called \textit{irreducible} if the only proper
subspace $H' \subseteq H$ with $GH' \subseteq H'$ is $\{0\}$. Two
$G$-spaces $H$ and $H'$ are called \textit{equivalent} if there is a
\textit{$G$-isometry} $\phi : H \to H'$, i.e. a linear isomorphism
with $\phi(af) = a\phi(f)$ for all $f \in H$ and $a \in G$ and
$(\phi(f), \phi(g)) = (f, g)$ for all $f, g \in H$.

By Maschke's theorem (cf. \cite[Theorem 2.4.1]{goodman-wallach-1998})
one can decompose $\C^X$ orthogonally into irreducible $G$-spaces:
\begin{equation}
\label{decomp cx}
\C^X = (H_{1, 1} \perp \ldots \perp H_{1, m_1}) \perp \ldots \perp
(H_{d, 1} \perp \ldots \perp H_{d, m_d}),
\end{equation}
where $H_{k,i}$ with $k = 1, \ldots, d$ and $i = 1, \ldots, m_k$ is an
irreducible $G$-space of dimension $h_k$ and where $H_{k,i}$ and
$H_{k',i'}$ are equivalent if and only if $k = k'$.

Let $\MA$ be the subalgebra of $\C^{X \times X}$ which is generated by
the permutation matrices $P_a \in \C^{X \times X}$ with $a \in G$
where
\begin{equation}
\label{Pa}
P_a(x,y) = 
\begin{cases}
1 & \text{if $a^{-1}x = y$,} \\ 
0 & \text{otherwise.}
\end{cases}
\end{equation}
Because of \eqref{decomp cx} the algebra $\MA$ decomposes as a complex
vector space in the following way
\begin{equation}
\label{decomp a}
\MA \cong \bigoplus_{k=1}^d \C^{h_k \times h_k} \otimes I_{m_k}.
\end{equation}

Recall that by $\MB$ we denote the matrices in $\C^{X \times X}$ which
are invariant under the group action of $G$. In other words, it is the
\textit{commutant} of $\MA$:
\begin{equation*}
\MB = \Comm(\MA) = \{B \in \C^{X \times X} : \mbox{$BA = AB$ for all
$A \in \MA$}\}.
\end{equation*}
The double commutant theorem \cite[Theorem
3.3.7]{goodman-wallach-1998} gives the following decomposition of
$\MB$ as a complex vector space:
\begin{equation}
\label{decomp b}
\MB \cong \bigoplus_{k=1}^d I_{h_k} \otimes \C^{m_k \times m_k}.
\end{equation}

Now we construct an explicit algebra isomorphism between the commutant
algebra $\MB$ and matrix algebra $\bigoplus_{k=1}^d \C^{m_k \times
  m_k}$. 

Let $e_{k, 1, l}$ with $l = 1, \ldots, h_k$ be an orthonormal basis of
the space $H_{k,1}$. Choose $G$-isometries $\phi_{k,i} : H_{k,1} \to
H_{k,i}$. Then, $e_{k,i,l} = \phi_{k,i}(e_{k,1,l})$ is an orthonormal
basis of $H_{k,i}$. Define the matrix $E_{k,i,j} \in \C^{X \times X}$
with $i, j = 1, \ldots, m_k$ by
\begin{equation*}
E_{k,i,j}(x, y) = \frac{1}{|X|} \sum_{l = 1}^{h_k} e_{k,i,l}(x) \overline{e_{k,j,l}(y)}.
\end{equation*}

The definition of these matrices depend on the choice of the
orthonormal basis, on the chosen $G$-isometries and on the chosen
decomposition \eqref{decomp cx}. The following proposition shows
the effect of different choices. 

\begin{proposition}
By $E_k(x,y)$ we denote the $m_k
\times m_k$ matrix $(E_{k,i,j}(x,y))_{i,j}$.

(a) The matrix entries $E_{k,i,j}(x, y)$ do not depend on the choice
of the orthonormal basis of $H_{k,1}$.

(b) The change of $\phi_{k,i}$ to $\alpha \phi_{k,i}$ with $\alpha \in
\C$, $|\alpha| = 1$,  simultaneously changes the $i$-th
row and $i$-th column in the matrix $E_k(x,y)$ by a multiplication
with $\alpha$ and $\overline{\alpha}$, respectively.

(c) The choice of another decomposition of $H_{k,1} \perp \ldots \perp
H_{k, m_k}$ as a sum of $m_k$ orthogonal, irreducible $G$-spaces
changes $E_k(x,y)$ to $UE_k(x,y)\overline{U}^t$ for some unitary matrix $U \in
\Uni(\C^{m_k})$.
\end{proposition}

\begin{proof}
This was proved in \cite[Theorem 3.1]{bachoc-vallentin-2006a} with the
only difference that there only the real case was considered. The
complex case follows mutatis mutandis.
\end{proof}

The following theorem shows that the map
\begin{equation}
\label{psi}
\varphi : \MB \to \bigoplus_{k=1}^d \C^{m_k \times m_k}
\end{equation}
mapping $E_{k,i,j}$ to the elementary matrix with the only non-zero
entry $1$ at position $(i,j)$ in the $k$-th summand $\C^{m_k \times
m_k}$ of the direct sum is an algebra isomorphism.

\begin{theorem}
The matrices $E_{k,i, j}$ form a basis of $\MB$ satisfying the equation
\begin{equation}
\label{multiplication}
E_{k,i,j} E_{k',i',j'} = \delta_{k,k'} \delta_{j,i'} E_{k,i,j'},
\end{equation}
where $\delta$ denotes Kronecker's delta.
\end{theorem}

\begin{proof}
The multiplication formula \eqref{multiplication} is a direct
 consequence of the orthonormality of the vectors $e_{k,i,l}$.  That
 $E_{k,i,j}$ is an element of $\MB$ follows from \cite[Theorem 3.1
 (c)]{bachoc-vallentin-2006a}. From \eqref{multiplication} it follows
 that the matrices $E_{k,i,j}$ are linearly independent, they span a
 vector space of dimension $\sum_{k = 1}^d m_k^2$. Hence, by
 \eqref{decomp b}, they form a basis of the commutant $\MB$.
\end{proof}

Now the expansion of the canonical basis $B_r$, with $r = 1, \ldots,
N$, in the basis $E_{k,i,j}$ with coefficients $p_r(k,i,j)$
\begin{equation}
\label{p functions}
B_r = \sum_{k = 1}^d \sum_{i,j = 1}^{m_k} p_r(k,i,j) E_{k,i,j}.
\end{equation}
yields
\begin{equation*}
\varphi(B_r) = \sum_{k = 1}^d \sum_{i,j = 1}^{m_k} p_r(k,i,j) \varphi(E_{k,i,j}).
\end{equation*}

\subsection{Orthogonality relation}

For the computation of the coefficients $p_r(k,i,j)$ the following
orthogonality relation is often helpful.

If we expand the basis $|X| E_{k,i,j}$ in the canonical basis $B_r$ we get a
relation which after normalization is inverse to \eqref{p functions}
\begin{equation}
\label{q functions}
|X| E_{k,i,j} = \sum_{r = 1}^N q_{k,i,j}(r) B_r.
\end{equation}

So we have an orthogonality relation between the $q_{k,i,j}$:

\begin{lemma}
\label{orel}
Let $v_r = |\{(x,y) \in X \times X : (x,y) \in R_r\}|$. Then,
\begin{equation}
\label{orth rel}
\sum_{r = 1}^N v_r q_{k,i,j}(r) \overline{q_{k',i',j'}(r)} =
\delta_{k,k'} \delta_{j,j'} \delta_{i,i'} |X|^2 h_k.
\end{equation}
\end{lemma}

\begin{proof}
Consider the sum $\sum_{x \in X} E_{k,i,j} E_{k',j',i'} (x,x)$. 

On the one hand it is equal to
\begin{equation*}
\sum_{x \in X} \delta_{k,k'} \delta_{j,j'} E_{k,i,i'} (x,x) =  \delta_{k,k'} \delta_{j,j'} \trace E_{k,i,i'},
\end{equation*}
and
\begin{equation*}
\trace E_{k,i,i'} = \sum_{l = 1}^{h_k} (e_{k,i,l}, e_{k,i',l}) = \delta_{i,i'} h_k,
\end{equation*}

On the other hand it is
\begin{equation*}
\sum_{x \in X} \sum_{y \in X} E_{k,i,j}(x,y) E_{k',j',i'}(y,x)
 =  \frac{1}{|X|^2} \sum_{r = 1}^N v_r q_{k,i,j}(r) \overline{q_{k',i',j'}(r)},
\end{equation*}
where we used the fact $E_{k',j',i'}(y,x) =
\overline{E_{k',i',j'}(x,y)}$ which follows from the
definition.
\end{proof}

The orthogonality relation gives a direct way to compute $p_r(k,i,j)$
once $q_{k,i,j}(r)$ is known: We have
\begin{equation}
\label{direct inverse q}
p_r(k,i,j) = \frac{v_r \overline{q_{k,i,j}(r)}}{|X| h_k},
\end{equation}
which follows by Lemma \ref{orel} and by \eqref{p functions} and
\eqref{q functions} because of
\begin{equation*}
\sum_{r = 1}^N p_r(k,i,j) q_{k',i',j'}(r) = |X| \delta_{k,k'} \delta_{i,i'} \delta_{j,j'}.
\end{equation*}

\subsection{Algorithmic issues}

We conclude this section by reviewing algorithmic issues for computing
$\varphi$. To calculate the isomorphism one has to perform the
following steps:
\begin{enumerate}
\item Compute the orthogonal decomposition \eqref{decomp cx} of $\C^X$
  into pairwise orthogonal, irreducible $G$-spaces $H_{k,i}$.
\item For every irreducible $G$-space $H_{k,1}$ determine an
orthonormal basis.
\item Find $G$-isometries $\phi_{k,i} : H_{k,1} \to H_{k,i}$.
\item Express the basis $B_r$ in the basis $E_{k,i,j}$.
\end{enumerate}

Only the first step requires an algorithm which is not classical. Here
one can use an algorithm of Babai and R\'onyai \cite{br-1990}. It is a
randomized algorithm running in expected polynomial time for computing
the orthogonal decomposition \eqref{decomp cx}. It requires the
permutation matrices $P_a$ given in \eqref{Pa} as input, where $a$
runs through a (favorably small) generating set of $G$. The other
steps can be carried out using Gram-Schmidt orthonormalization and
solving systems of linear equations.

\section{Block diagonalization of the Terwilliger algebra}
\label{terwilliger}

The symmetric group $S_n$ acts on the set $X = \{0,1\}^n$ of binary
vectors with length $n$ by $\sigma(x_1, \ldots, x_n) = (x_{\sigma(1)},
\ldots, x_{\sigma(n)})$, i.e.\ by permuting coordinates. In
\cite{schrijver-2005} Schrijver determined the block diagonalization
of the algebra $\MB$ of $X \times X$-matrices invariant under this
group action. The algebra $\MB$ is called the \textit{Terwilliger
  algebra of the binary Hamming scheme}. Now we shall derive a block
diagonalization in the framework of the previous section. In this case it
is possible to work over the real numbers only because all irreducible
representations of the symmetric group are real.

Under the group action the set $X$ splits into $n + 1$ orbits $X_0,
\ldots, X_n$ where $X_m$ contains the elements of $\{0,1\}^n$ having
Hamming weight $m$, i.e.\ elements which one can get from the binary
vector $1^m0^{n-m}$ by permuting coordinates. So we have the
orthogonal decomposition of the $S_n$-space $\R^X$ into
\begin{equation*}
\R^{X} = \R^{X_0} \perp \ldots \perp \R^{X_n}.
\end{equation*}

It is a classical fact (cf. \cite[Theorem 2.10]{dunkl-1976}) that the
$S_n$-space $\R^{X_m}$ decomposes further into
\begin{equation*}
\R^{X_m} = 
\left\{
\begin{array}{ll}
  H_{0,m} \perp \ldots \perp H_{m, m}, & \mbox{when $0 \leq m \leq \lfloor n/2 \rfloor$,}\\
H_{0,m} \perp \ldots \perp H_{n-m, m}, & \mbox{otherwise.}
\end{array}
\right.
\end{equation*}
where $H_{k,m}$ are irreducible $S_n$-spaces which correspond to the
irreducible representation of $S_n$ given by the partition $(n-k, k)$
(cf. \cite[Chapter 2]{sagan-2001}). Its dimension is $h_k =
\binom{n}{k} - \binom{n}{k-1}$.

Thus, the matrices $E_{k,i,j}$, with $k = 0, \ldots, \lfloor n/2
\rfloor$, which correspond to the isotypic component $H_{k,k} \perp
\ldots \perp H_{k,n-k}$ of $\R^X$ of type $(n-k,k)$ are conveniently
indexed by $i,j = k, \ldots, n-k$. Since $E_{k,j,i}$ is the transpose
of $E_{k,i,j}$ we only need to consider the case $k \leq i \leq j \leq n - k$.

To determine $E_{k,i,j}(x,y)$ we rely on the papers \cite{dunkl-1976}
and \cite{dunkl-1979} of Dunkl. We recall the facts and notation which
we will need from them. Let $T_k : S_n \to \Ort(\R^{h_k})$ be an
orthogonal, irreducible representation of $S_n$ given by the partition
$(n-k, k)$. By $H, K$ we denote the subgroups $H = S_j \times S_{n-j}$
and $K = S_i \times S_{n-i}$ of $S_n$.  Let $V_k \subseteq \R^{S_n}$
be the vector space spanned by the function $(T_k)_{rs}$, with $1 \leq
r, s \leq h_k$, which are the matrix entries of $T_k$:
$(T_k)_{rs}(\pi) = [T_k(\pi)]_{rs}$. A function $f \in V_k$ is called
\textit{$H$-$K$-invariant} if $f(\sigma \pi \tau) = f(\pi)$ for all
$\sigma \in H$, $\pi \in S_n$, $\tau \in K$. In \cite[\S
4]{dunkl-1976} and \cite[\S 4]{dunkl-1979} Dunkl computed the
$H$-$K$-invariant functions of $V_k$. These are all real multiples of
\begin{equation*}
\psi_{k, H-K}(\pi) = \frac{(-j)_k(i-n)_k}{(-i)_k(j-n)_k} Q_k(v(\pi); -(n-i)-1, -i-1, j),
\end{equation*}
where $(a)_0 = 1$, $(a)_k = a(a+1)\ldots(a+k-1)$, and where,
\begin{equation*}
Q_k(x; -a-1, -b-1, m) = \frac{1}{\binom{m}{k}} \sum_{j=0}^k (-1)^j \frac{\binom{b-k+j}{j}}{\binom{a}{j}} \binom{m-x}{k-j} \binom{x}{j},
\end{equation*}
are \textit{Hahn polynomials} (for integers $m,a,b$ with $a \geq m, b \geq m
\geq 0)$, and where
\begin{equation*}
v(\pi) = i-|\pi\{1,\ldots,i\} \cap \{1, \ldots, j\}|.
\end{equation*}
The polynomials $Q_k(x) = Q_k(x; -a-1, -b-1, m)$ are the orthogonal
polynomials for the weight function $\binom{a}{x}\binom{b}{m-x}$, $x =
0, 1, \ldots, m$, normalized by $Q_k(0) = 1$. For more information
about Hahn polynomials we refer to \cite{km-1961}.

We will need the square of the norm of $\psi_{k, H-K}$ which is given
in \cite[before Proposition 2.7]{dunkl-1979}:
\begin{equation*}
(\psi_{k, H-K}, \psi_{k, H-K}) = \frac{\psi_{k, H-K}(\mathrm{id})}{h_k} =
\frac{(-j)_k(i-n)_k}{(-i)_k(j-n)_k h_k}.
\end{equation*}

Let $e_{k,i,1}, \ldots, e_{k,i,h_k}$ be an orthonormal basis of
$H_{k,i}$. We get an orthogonal, irreducible representation
$T_{k,i} : S_n \to \Ort(\R^{h_k})$ by
\begin{equation*}
\pi (e_{k,i,l}) = \sum_{l'=1}^{h_k} [T_{k,i}(\pi)]_{l',l} e_{k,i,l'}.
\end{equation*}
Consider the function
\begin{equation*}
z_{k, i, j}(\pi) = E_{k,i,j}(\pi(1^i0^{n-i}), 1^j0^{n-j}).
\end{equation*}
This is an $H$-$K$-invariant function because $E_{k,i,j} \in \MB$.  It
lies in $V_k$ because vector spaces spanned by matrix entries of two
equivalent irreducible representations coincide. Thus, $z_{k, i, j}$
is a real multiple of $\psi_{k, H-K}$. By computing the squared norm
of $z_{k, i, j}$ we determine this multiple up to sign:
\begin{eqnarray*}
(z_{k, i, j}, z_{k, i, j}) & = & \frac{1}{n!} \sum_{\pi \in S_n} z_{k, i, j}(\pi) z_{k, i, j}(\pi)\\
& = & \frac{1}{\binom{n}{i} 2^n} \sum_{l = 1}^{h_k} (e_{k, i, l}(1^j 0^{n-j}))^2\\
& = & \frac{1}{\binom{n}{i}} E_{k,j,j}(1^j0^{n-j}, 1^j 0^{n-j}).
\end{eqnarray*}
Here we used that $e_{k,i,l}$ is an orthonormal basis of $H_{k,i}$
where the inner product is $(f,g) = \frac{1}{2^n} \sum_{x \in X_i}
f(x) g(x)$.

All diagonal entries belonging to $X_j \times X_j$ of $E_{k,j,j}$
coincide and all others are zero, so
$\binom{n}{j}E_{k,j,j}(1^j0^{n-j}, 1^j 0^{n-j})$ is the trace of
$E_{k,j,j}$ which equals its rank $h_k$. Hence, $(z_{k, i, j}, z_{k,
i, j}) = h_k (\binom{n}{i}\binom{n}{j})^{-1}$. So we have determined
$E_{k,i,j}$ up to sign. To adjust the signs it is enough to ensure
that the multiplication formula \eqref{multiplication} is satisfied.

So putting it together, we have proved the following theorem.

\begin{theorem}
\label{Ekst}
For $x, y \in X$ define $v(x,y) = |\{l \in \{1, \ldots, n\} : x_l =
1, y_l = 0\}|$.  For $k = 0, \ldots, \lfloor n/2 \rfloor$ and $i, j
= k, \ldots, n-k$ with $i \leq j$ we have
\begin{equation*}
\begin{split}
E_{k,i,j}(x,y) & = \frac{h_k}{(\binom{n}{i}\binom{n}{j})^{1/2}}
\left(\frac{(-j)_k(i-n)_k}{(-i)_k(j-n)_k}\right)^{-\frac{1}{2}} \cdot\\
& \quad\qquad Q_k(v(x,y); -(n-i)-1, -i-1, j),
\end{split}
\end{equation*}
when $x \in X_i$, $y \in X_j$. In the case $x \not\in X_i$ or $y
\not\in X_j$ we have $E_{k,i,j}(x,y) = 0$. Furthermore, $E_{k,j,i} =
(E_{k,i,j})^t$.
\end{theorem}

Finally, to find the desired algebra isomorphism \eqref{isomorphism}
we determine the values of $p_r(k,i,j)$ by formula \eqref{direct
inverse q}. We represent the orbits $R_1, \ldots, R_N$ by triples
$(r,s,d)$: Two pairs $(x,y), (x',y') \in X \times X$ are equivalent
whenever $x, x' \in X_r$, $y, y' \in X_s$, and $v(x,y) = v(x',y') =
d$. Then,
\begin{equation*}
p_{r,s,d}(k,i,j) = \frac{v_{r,s,d}E_{k,i,j}(x,y)}{h_k},
\end{equation*}
where
\begin{equation*}
v_{r,s,d} = \binom{n}{d} \binom{n-d}{r-d} \binom{n-r}{s-s+d}.
\end{equation*}

\begin{remark}
In a similar way one can give an interpretation of the block
diagonalization of the Terwilliger algebra for nonbinary codes which
was computed in \cite{gst-2006}. Using \cite[Theorem 4.2]{dunkl-1976}
one can show the matrix entries are, up to scaling factors, products
of Hahn polynomials and Krawtchouk polynomials.
\end{remark}

\section*{Acknowledgements}

We thank Christine Bachoc and the anonymous referee for very helpful
remarks and suggestions.

\end{document}